\documentclass[12pt,draft]{amsart}
\usepackage[all]{xy}
\usepackage{amsfonts}
\usepackage{amssymb}

\newtheorem{teo}{Theorem}[section]
\newtheorem{lem}[teo]{Lemma}

\theoremstyle{definition}
\newtheorem{dfn}[teo]{Definition}

\def\<{\langle}
\def\>{\rangle}
\def\ss{\subset}

\def\a{\alpha}

\def\e{\varepsilon}

\def\r{\rho}

\def\f{{\varphi}}

\def\C{{\mathbb C}}

\def\Z{{\mathbb Z}}

\def\Id{\operatorname{Id}}

\def\Prim{\operatorname{Prim}}

\def\Tr{\operatorname{Tr}}

\newcommand{\ov}[1]{\overline{#1}}

\newcommand{\wh}[1]{\widehat{#1}}

\newcommand{\Mat}[4]{\left( \begin{array}{cc}
                            #1 & #2 \\
                            #3 & #4
                      \end{array} \right)}
\newcommand{\Matt}[9]{\left( \begin{array}{ccc}
                            #1 & #2 & #3 \\
                            #4 & #5 & #6 \\
                            #7 & #8 & #9
                      \end{array} \right)}

\def\Fix{\operatorname{Fix}}

\def\N{{\mathbb N}}
\def\T{{\mathbb T}}

\begin{document}

\title[An example of twisted Burnside theorem]
{Twisted Burnside theorem for type II${}_1$ groups:
an example}

\author[A.~Fel'shtyn]{Alexander Fel'shtyn}
\address{Instytut Matematyki, Uniwersytet Szczecinski,
ul. Wielkopolska 15,
70-451 Szczecin, Poland
and Department of Mathematics, Boise State University,
1910 University Drive, Boise, Idaho, 83725-155, USA }
\email{felshtyn@diamond.boisestate.edu, felshtyn@mpim-bonn.mpg.de}
\urladdr{
http://www.mat.univ.szczecin.pl/en/staff/felshtyn\_a.html}
\author[E.~Troitsky]{Evgenij Troitsky}
\thanks{This a joint research completed during the stay
at Max-Planck-Institut f\"ur Mathematik (Bonn) in Spring 2004
in partial relation to the activity on Algebraic and Topological
Dynamics.
The authors are grateful to MPIM for support and hospitality.
\newline
The second author is partially supported by
RFBR Grant 05-01-00293 and Grant PH$\Pi$.2.1.1.5055.
\newline
The third author is partially supported by
RFBR Grant 05-01-00899 (Russia) and CRDF Grant RUM-2622.ST.04 (US).
\newline
The authors are indebted to the referee for very useful suggestions.}
\address{Dept. of Mech. and Math., Moscow State University,
119992 GSP-2  Moscow, Russia}
\email{troitsky@mech.math.msu.su}
\urladdr{
http://mech.math.msu.su/\~{}troitsky}
\author[A.~Vershik]{Anatoly Vershik}
\address{St.Petersburg Department of Steklov Institute of Mathematics,
27 Fontanka, St. Petersburg 191011, Russia}
\email{vershik@pdmi.ras.ru}
\urladdr{http://www.pdmi.ras.ru/\~{}vershik}

\begin{abstract}
The purpose of the present paper is to discuss
the following conjecture of Fel'shtyn and Hill, which is a
generalization of the classical Burnside theorem:

Let $G$ be a countable discrete group, $\phi$ its  automorphism,
$R(\phi)$ the number of $\phi$-conjugacy classes
(Reidemeister number),
$S(\phi)=\# \Fix (\wh\phi)$ the
number of $\phi$-invariant equivalence classes of irreducible
unitary representations. If one of $R(\phi)$ and
$S(\phi)$ is finite, then it is equal to the other.

This conjecture plays a very important role in the theory of
twisted conjugacy classes
having a long history (see \cite{Jiang}, \cite{FelshB})
and has very serious consequences
in Dynamics, while its proof needs rather fine results from
Functional and Non-commutative Harmonic Analysis.
It was proved for finitely generated groups of type I in \cite{FelTro}.

In the present paper this conjecture is disproved for non-type I groups.
More precisely, an example of a group and its automorphism is
constructed such that the number of fixed irreducible
representations is greater than the Reidemeister number. But the
number of fixed finite-dimensional representations
(i.e. the number of invariant finite-dimensional characters)
in this example
coincides with the Reidemeister number.

The directions for search of an appropriate formulation are indicated
(another definition of the dual object).
\end{abstract}

\maketitle
\tableofcontents

\section{Introduction: motivation, history, and current state}

\begin{dfn}
Let $G$ be a countable discrete group and $\phi: G\rightarrow G$ an
endomorphism.
Two elements $x,x'\in G$ are said to be
 $\phi$-{\em conjugate} or {\em twisted conjugate,}
iff there exists $g \in G$ with
$$
x'=g  x   \phi(g^{-1}).
$$
We shall write $\{x\}_\phi$ for the $\phi$-{\em conjugacy} or
{\em twisted conjugacy} class
 of the element $x\in G$.
The number of $\phi$-conjugacy classes is called the {\em Reidemeister number}
of an  endomorphism $\phi$ and is  denoted by $R(\phi)$.
If $\phi$ is the identity map then the $\phi$-conjugacy classes are the usual
conjugacy classes in the group $G$.
\end{dfn}

If $G$ is a finite group, then the classical Burnside theorem (see e.g.
\cite[p.~140]{Kirillov})
says that the number of
classes of irreducible representations is equal to the number of conjugacy
classes of elements of $G$.  Let $\wh G$ be the {\em unitary dual} of $G$,
i.e. the set of equivalence classes of unitary irreducible
representations of $G$.

Therefore, by the Burnside's theorem, if $\phi$ is the identity automorphism
of any finite group $G$, then we have
 $R(\phi)=\#\Fix(\wh\phi)$, where $\wh\phi[\rho]=[\rho\circ\phi]$.

One of the main achievements in the field till now is the following result.

\begin{teo}[\cite{FelTro}]
\label{teo:mainth1} Let $G$ be a finitely generated discrete group
of type {\rm I}, $\phi$ one of its endomorphism, $R(\phi)$ the
number of $\phi$-conjugacy classes, and $S(\phi)=\# \Fix
(\wh\phi)$ the number of $\wh\phi$-invariant equivalence classes
of irreducible unitary representations. If one of $R(\phi)$ and
$S(\phi)$ is finite, then it is equal to the other.
\end{teo}

The research is motivated not only by a natural desire to
extend the classical Burnside theorem to the case of infinite
groups and twisted conjugacy classes, but also by dynamical
applications. Namely, the identification of Reidemeister number
with a number of fixed points in a natural way, has some very
interesting consequences.

More precisely,
let $\mu(d)$, $d\in\N$, be the {\em M\"obius function},
i.e.
$$
\mu(d) =
\left\{
\begin{array}{ll}
1 & {\rm if}\ d=1,  \\
(-1)^k & {\rm if}\ d\ {\rm is\ a\ product\ of}\ k\ {\rm distinct\ primes,}\\
0 & {\rm if}\ d\ {\rm is\ not\ square-free.}
\end{array}
\right.
$$

\begin{teo}[Congruences for the Reidemeister numbers \cite{FelTro}]\label{teo:mainth3}
Let $\phi:G$ $\to G$ be an endomorphism of a countable discrete group $G$
such that
all numbers $R(\phi^n)$ are finite and let $H$ be a subgroup
 of $G$ with the properties
$$
  \phi(H) \subset H,
$$
$$
  \forall x\in G \; \exists n\in \N \hbox{ such that } \phi^n(x)\in H.
$$
If for the pair  $(H,\phi^n)$ the twisted Burnside theorem holds,
i.e. this pair satisfies the conclusion of
Theorem~{\rm~\ref{teo:mainth1}}
for any $n\in\N$,
then one has for all $n$,
 $$
 \sum_{d\mid n} \mu(d)\cdot R(\phi^{n/d}) \equiv 0 \mod n.
 $$
\end{teo}

These theorems were proved previously in a  special case of
a direct sum of an
Abelian finitely generated group and finite group \cite{FelHill,FelBanach}.
Remember that the difference between result of \cite{FelHill,FelBanach}
and Theorem \ref{teo:mainth1} is not so big because by theorem of E.Thoma a
countable finitely generated group $G$ is a group of type I iff
it is a finite extension of an Abelian group.
The conjecture under discussion was formulated in \cite{FelHill}.
We refer to \cite{FelshB,FelTro} for more detail and topological
applications.

On the other hand, one can introduce the number
$R_*(\phi)$ of "Reidemeister classes related to
twisted invariant functions on $G$ from the Fourier-Stieltjes
algebra $B(G)$", or more precisely, the dimension of the space
of twisted invariant functions on $G$ which can be extended up to
bounded functionals on the group algebra $C^*(G)$. Let $S_*(\phi)$
be the sum of codimensions of subspaces  $L_I\ss C^*(G)/I$, where
$L_I$ is generated by elements of the form $a-L_g a L_{\phi(g^{-1})}$
and $I$ runs over the \emph{Glimm spectrum} $T$ of $G$.
Let us remind this notion.
Let $Z$ be the center of a C*-algebra $A$ and $
\wh Z$ its space of maximal ideals
equipped with the standard topology. If  $I\in {\mathcal P}:=\Prim (A)$
(the space of kernels of unitary irreducible representations), then
$Z\cap I \in \wh Z$ (this follows from the fact that
the restriction onto $Z$ of an irreducible representation with kernel
$I$ gives rise to a homomorphism
$Z \to \C$, and hence is a maximal ideal $Z\cap I$). We obtain a map
$f: {\mathcal P}\to \wh Z$. Suppose $T:=f({\mathcal P})$. For each $x\in T$ consider the ideal
$I_x:=\cap I$, $f(I)=x$, (\emph{Glimm ideal}) and the field of algebras
$A/I_x$. We have the map $a\mapsto \{x\mapsto a+I_x\}$ from the algebra
$A$ to the algebra of sections of the mentioned field. An important
result of
\cite{DaunsHofmann} asserts that this map is an isomorphism.
The map
$f:{\mathcal P}\to T$ is universal with respect to continuous maps  $g:{\mathcal P}\to S$
to Hausdorff spaces, i.e. any such map can be represented under the
form
$h\circ f$ for some continuous $h:T\to S$.
The space $T$ is compact for a unital algebra.

We call $S_*(\phi)$ the
number of generalized fixed points of $\wh\phi$ on the
Glimm spectrum of $G$.

 \begin{teo}[{weak twisted Burnside theorem, \cite{ncrmkwb}}]\label{teo:weakburn}
The number $R_*(\phi)$  is  equal to the number $S_*(\phi)$ of generalized
fixed points of $\wh\phi$ on the Glimm spectrum of $G$,
if one of $R_*(\phi)$ and
$S_*(\phi)$ is finite.
 \end{teo}
 This result
allows to obtain the strong form of twisted
Burnside theorem $R(\phi)=S(\phi)$ in a number of cases.

\medskip
The interest in twisted conjugacy relations has its origins, in particular,
in the Nielsen-Reidemeister fixed point theory (see, e.g. \cite{Jiang,FelshB}),
in Selberg theory (see, eg. \cite{Shokra,Arthur}),
and  Algebraic Geometry (see, e.g. \cite{Groth}).

The congruences give some necessary conditions for the realization problem
for Reidemeister numbers in topological dynamics.

Let us remark that it is known that the Reidemeister number of an
endomorphism of a finitely generated Abelian group is
finite iff $1$ is not in the
spectrum of the restriction of this endomorphism to the free part of the group
(see, e.g. \cite{Jiang}).  The Reidemeister number
is infinite for any automorphism of a  non-elementary
Gromov hyperbolic group \cite{FelPOMI,ll}
as well as for any
injective endomorphism of Baumslag-Solitar group \cite{FelGon} (see also \cite{FelTroObzo}).

The main result of the present paper is the following statement.

\begin{teo}
There exists an amenable solvable group $G$ not of type I and its
automorphism $\phi$ such that its Reidemeister number $R(\phi)$
is finite but does not coincide with the number $S(\phi)$ of
fixed points of $\wh\phi$ on $\wh G$.
\end{teo}

This example is very important for the further attack onto
the problem, because this group is "situated between" the groups
of type I, for which the conjecture is true \cite{FelTro},
and Gromov hyperbolic
groups, for which Reidemeister numbers are always infinite.

The ideas for the further study arising from the example from
this paper and from \cite{ncrmkwb} are discussed in the last
section below.

\section{Some technical preliminaries}

Let $G$ be a semidirect product of $\Z^2$ and $\Z$ by Anosov automorphism
$\a$ with the matrix $A=\Mat 2111$ (of course, our results remain valid for
any hyperbolic element of $SL(2,\Z)$ to use as $A$).
It consists by the definition of triples
$((m,k),n)$ of integers with the following multiplication low:
$$
((m,k),n)*((m',k'),n')=((m,k)+\a^n(m',k'),n+n').
$$
In particular,
$$
((m,k),0)*((0,0),n)=((m,k),n).
$$
Elements of  $G$ can be written also as matrices
$$
\bmatrix \alpha^n & x\\ 0 & 1\endbmatrix, \qquad n\in\Z, \quad x=(m,k)\in\Z^2.
$$
The inverse of $((m,k),n)$ is $(-\a^{-n}(m,k),-n)$. Indeed,
$$
((m,k),n)*(-\a^{-n}(m,k),-n)=((m,k)-\a^n \a^{-n}(m,k),n-n)=((0,0),0).
$$

The group $G$ is a solvable (hence, amenable) group which is not
of type I. Its regular representation is factorial.

Let us define an automorphism $\phi:G\to G$ by
$$
\phi((m,k),n)=((k,-m),-n),
$$
i.e. the action on $\Z^2$ is defined by automorphism $\mu$ with the matrix
$M=\Mat 01{-1}0$, and on $\Z$ by $n\mapsto -n$. The map $\phi$ is clearly
a bijection,
\begin{multline*}
\phi(((m,k),n)*((m',k'),n')=\phi((m,k)+\a^n(m',k'),n+n'))\\
=
((k,-m)+\mu\a^n(m',k'),-n-n'),
\end{multline*}
\begin{multline*}
\phi((m,k),n)*\phi(((m',k'),n')=((k,-m),-n)*((k',-m'),-n') \\
=((k,-m)+\a^{-n}(k',-m'),-n-n').
\end{multline*}
Hence, to prove that $\phi$ is an automorphism, we need $\mu\a^n=\a^{-n}\mu$.
This follows from $\mu\a=\a^{-1}\mu$. The further results in this direction
can be found in \cite{gowon}.

       One can hope to find $\wh\phi$-fixed
       irreducible representations (at least finite-dimensional
       ones) from those orbits of $\alpha$ on the dual
       torus $\T^2$ to the normal subgroup $\Z^2$, which are invariant
       under the dual action of $\mu$, using appropriate
       cocycles. This will be done in the last section. We will
       find four such representations and we will explain why
       that is all. Some motivation for this way of construction
       can be found in \cite[Ch.~17, \S~1]{baron}.

\section{Description of Reidemeister classes}

Let us find the Reidemeister classes of $\phi$, i.e. the classes of the
equivalence relation $h\sim gh\phi(g^{-1})$. We will prove the
following statement.

\begin{lem}
For the above $G$ and $\phi$ one has $R(\phi)=4$.
\end{lem}

Also we will obtain a description of these four Reidemeister classes.

For $h=((m,k),n)$ and
$g=((x,y),z)$ the right hand side of the relation takes the following
form:
$$
((x,y)+ \a^z(m,k),z+n)*(-\mu\a^{-z}(x,y),z)=
$$
$$
((x,y)+ \a^z(m,k)-\a^{z+n}\mu\a^{-z}(x,y),2z+n)=
$$
\begin{equation}\label{eq:formconj}
=(\a^z\{(m,k) + (\Id - \a^n\mu) \a^{-z}(x,y)\}, 2z+n).
\end{equation}
Let us call {\em level} $n$ (of $G$) the coset $L_n$ of $\Z^2\ss G$
of all elements of the form $((m,k),n)$. Let us first take an element
$((m,k),0)$ from the level $0$ and describe elements from the same
level, being equivalent to it. By (\ref{eq:formconj}) in this case
$z=n=0$ and they have the form
$$
((m,k) + (\Id - \mu)(x,y), 0)=((m + (x-y),k + (x+y)),0),
$$
where $\Id - \mu$ has the matrix $\Mat 1{-1}11$. Hence, the level $0$
has intersections with 2 Reidemeister classes, say, $B_1$ and $B_2$.
The first intersection $B_1\cap L_0$
is formed by elements $((u,v),0)$ with odd $u+v$,
and $B_2\cap L_0$ --- with even $u+v$. The elements from the other
levels, which are equivalent to $((m,k),0)$, have the form
\begin{equation}\label{eq:even1}
(\a^z\{(m,k) + (\Id - \mu) \a^{-z}(x,y)\}, 2z).
\end{equation}
This means that $B_1$ and $B_2$ enter only even levels.
Also, since $\a$ is an automorphism, we can rewrite (\ref{eq:even1})
as
\begin{equation}\label{eq:even2}
(\a^z\{(m,k) + (\Id - \mu)(u,v)\}, 2z).
\end{equation}
with arbitrary integers $u$ and $v$. This means, that the intersections
$B_i\cap L_{2z}$ have the form $\a^z(B_i)$, $i=1,2$. In particular,
the other Reidemeister classes do not enter even levels.

In a similar way, the elements of $L_1$ equivalent to $((m,k),1)$
have the form
$$
((m,k) + (\Id - \a\mu)(x,y), 1)=((m + (2x-2y),k + x),1).
$$
This means, that $L_1$ enters 2 classes: $B_3$ is formed by elements
with even first coordinate, and $B_4$ --- with odd one.
The elements from the other
levels, which are equivalent to $((m,k),1)$, have the form
\begin{equation}\label{eq:odd1}
(\a^z\{(m,k) + (\Id - \a\mu) \a^{-z}(x,y)\}, 2z+1).
\end{equation}
Since $\a$ is an automorphism, we can rewrite (\ref{eq:odd1})
as
\begin{equation}\label{eq:odd2}
(\a^z\{(m,k) + (\Id - \a \mu)(u,v)\}, 2z+1).
\end{equation}
with arbitrary integers $u$ and $v$. This means, that the intersections
$B_i\cap L_{2z+1}$ have the form $\a^z(B_i)$, $i=3,4$. In particular,
these four classes cover $G$ and $R(\phi)=4$.

To obtain a complete description of $B_i$ let us remark that directly
from the definition of $\a$
$$
\a(x,y)=(2x+y,x+y)
$$
one has the following properties:

\begin{itemize}
\item $\a$ maps the set of elements with an even
(resp., odd) sum of coordinates onto the set of elements with an even
(resp., odd) second coordinate,

\item $\a$ maps the set of elements   with an even
(resp., odd) second coordinate onto the set of elements with an even
(resp., odd) first coordinate,

\item $\a$ maps the set of elements with an even
(resp., odd) first coordinate
 onto the set of elements with an even
(resp., odd) sum of coordinates.
\end{itemize}

Hence, the elements $((m,k),n)$ in intersections $B_i\cap L_j$
are of the form

\begin{center}
\begin{tabular}{|c|c|c|c|c|}
\hline
$i$  & 1 & 2 & 3 & 4 \\
\hline
\hline
$j\equiv 0 \mod 6$ & $m+k$ is even &$m+k$ is odd & $\emptyset$ & $\emptyset$ \\
\hline
$j\equiv 1 \mod 6$ &$\emptyset$ &$\emptyset$ & $m$ is even & $m$ is odd \\
\hline
$j\equiv 2 \mod 6$ & $k$ is even & $k$ is odd & $\emptyset$ & $\emptyset$ \\
\hline
$j\equiv 3 \mod 6$ &$\emptyset$ &$\emptyset$  & $m+k$ is even & $m+k$ is odd \\
\hline
$j\equiv 4 \mod 6$ & $m$ is even & $m$ is odd & $\emptyset$ & $\emptyset$ \\
\hline
$j\equiv 5 \mod 6$ &$\emptyset$ &$\emptyset$ & $k$ is even & $k$ is odd \\
\hline
\end{tabular}
\end{center}

\section{Fixed points and twisted invariant functionals}

Now we want to study the fixed points of the homeomorphism
$\wh\phi:\wh G\to \wh G$, $[\r]\mapsto [\r\phi]$ of the unitary dual.
We will prove the following statement.

\begin{lem}
There are exactly four finite-dimensional representations, which are
fixed points of $\wh\phi$. There is at least one
infinite-dimensional representation, which is a
fixed point of $\wh\phi$.
\end{lem}

Let us start from the finite-dimensional representations. As it was
shown in \cite{FelTro} in this case there exists exactly one twisted-invariant
functional on $L^1(G)$, or $\phi$-central $L^\infty$ function, coming
from a twisted-invariant functional on $\r(L^1(G))\cong M(\dim \r,\C)$
(up to scaling), namely
\begin{equation}\label{eq:trsg}
\f_\r:g\mapsto \Tr(S\r(g)),
\end{equation}
where $S$ is the intertwining operator between $\r$ and $\r\phi$.
This is an appropriate notion of character of a $\phi$-invariant
representation.

First, we have to find $\mu$-invariant finite $\a$-orbits on $\T^2$.
One can notice that
$$
\det(A^n-M)=\det\!\Mat a{b-1}{b+1}d\!=
ad-b^2+1=
\det A^n + 1=2,
$$
where  $A^n\!=\!\Mat abbd$,
for any $n$.
Hence, the mentioned orbits are formed by points with
coordinates $0$ and $1/2$. We have $2$ orbits: one of them consists of
$1$ point $(0,0)$ and gives rise to $1$-dimensional trivial representation
$\r_1$,
and the other consists of $A_1=(0,1/2)$, $A_2=(1/2,0)$ and
$A_3=(1/2,1/2)$ and
gives rise to a $3$-dimensional (irreducible) representation $\r_2$.
Also, one has the following $1$-dimensional representation $\pi$:
$$
\pi((m,k),2n)=1,\qquad \pi((m,k),2n+1)=-1.
$$
So, we have $4$ representations
$$
\r_1,\quad \r_2, \quad \pi,\quad \r_2 \otimes \pi.
$$
These representations are irreducible finite dimensional (for the last one
this follows from the evident observation $(\r_2 \otimes \pi)(g)=\pm \r_2(g)$,
hence the space of commuting operators is still isomorphic to $\C$).
We claim that these representations
give rise via (\ref{eq:trsg}) to $4$ linear independent
twisted-invariant functionals. In particular, there is no more
finite-dimensional fixed points of $\wh\phi$. More precisely, this
follows from the general theory of \cite{FelTro}, since these finite-dimensional
representations are well separated. But we will calculate directly for the
self-containing of the presentation here, and also because we need the explicit
form of these functionals.

Clearly,
\begin{equation}\label{eq:trsg1}
\f_{\r_1}\equiv 1, \quad \f_\pi =\left\{
\begin{array}{rl}
1, & \mbox{ on } \cup_n L_{2n}=B_1\cup B_2,\\
-1,& \mbox{ on } \cup_nL_{2n+1}=B_3\cup B_4,
\end{array}
\right.
\quad
\f_{\r_2\otimes \pi}=\f_{\r_2}\cdot \f_\pi.
\end{equation}
Let us find $\f_{\r_2}$. In the space $L^2(A_1,A_2,A_3)$ we take the base
$\e_1$, $\e_2$, $\e_3$ of characteristic functions in these points.
One has
$$
\a(\e_1)=\e_3,\qquad \a(\e_3)=\e_2,\qquad \a(\e_2)=\e_1,
$$
$$
\mu(\e_1)=\e_2,\qquad \mu(\e_2)=\e_1, \qquad \mu(\e_3)=\e_3.
$$
The representation (see \cite[Ch.~17, \S~1]{baron}) is defined by:
$$
\r_2(m,k,0)(\e_i)=\chi_{A_i}(m,k)\cdot \e_i,\qquad
\r_2(0,0,n)(\e_i)= \a^{-n} (\e_i)=\e_{i+n \mod 3},
$$
$$
\r_2(m,k,0)(\e_1)=e^{2\pi i(0\cdot m+ 1/2\cdot k)}\cdot \e_1
= e^{\pi i k} \e_1,
$$
$$
\r_2(m,k,0)(\e_2)=e^{2\pi i(1/2 \cdot m+ 0\cdot k)}\cdot \e_2
= e^{\pi i m} \e_2,
$$
$$
\r_2(m,k,0)(\e_3)=e^{2\pi i(1/2 \cdot m+ 1/2\cdot k)}\cdot \e_3
= e^{\pi i (m+k)} \e_3.
$$
The representation $\wh\phi\r_2$ is defined by
$$
\wh\phi\r_2(0,0,n)(\e_i)=\r_2(0,0,-n)(\e_i)
=\e_{i-n \mod 3},
$$
$$
\wh\phi\r_2(m,k,0)(\e_1)=\r_2(k,-m,0)(\e_1)
= e^{- \pi i m} \e_1= e^{\pi i m} \e_1,
$$
$$
\wh\phi\r_2(m,k,0)(\e_2)=\r_2(k,-m,0)(\e_2)
= e^{\pi i k} \e_2,
$$
$$
\wh\phi
\r_2(m,k,0)(\e_3)=\r_2(k,-m,0)(\e_3)
= e^{\pi i (k-m)} \e_3= e^{\pi i (m+k)} \e_3.
$$
The intertwining operator is induced by $\phi$ and has
the matrix $S=\Matt 010100001$. Hence,
$$
\f_{\r_2}(m,k,n)=\Tr(S\r_2(m,k,0)\r_2(0,0,n))=
$$
$$
=
\Tr\left[ \Matt 010100001
\Matt {e^{\pi i k}}000{e^{\pi i m}}000{e^{\pi i(k+m)}}
\Matt 001100010 ^n
\right]=
$$
$$
=
\Tr\left[
\Matt 0{e^{\pi i k}}0{e^{\pi i k}}0000{e^{\pi i(k+m)}}
{\Matt 001100010}^n
\right].
$$
For $n\equiv 0\mod 3$
\begin{eqnarray*}
\f_{\r_2}(m,k,n)&=&\Tr\left[
\Matt 0{e^{\pi i m}}0{e^{\pi i k}}0000{e^{\pi i(k+m)}}
\right]=e^{\pi i(k+m)}\\
&=&\left\{
\begin{array}{rl}
1,& \mbox{ if } m+k \mbox{ is even }\\
-1,& \mbox{ if } m+k \mbox{ is odd }
\end{array}
\right. ,
\end{eqnarray*}
for $n\equiv 1\mod 3$
$$
\f_{\r_2}(m,k,n)=\Tr\left[
\Matt 0{e^{\pi i m}}0{e^{\pi i k}}0000{e^{\pi i(k+m)}}
\Matt 001100010
\right]=
$$
$$
=
\Tr\Matt {e^{\pi i m}}0000{e^{\pi i k}}0{e^{\pi i(k+m)}}0
=e^{\pi i m}=\left\{
\begin{array}{rl}
1,& \mbox{ if } m \mbox{ is even }\\
-1,& \mbox{ if } m \mbox{ is odd }
\end{array}
\right. ,
$$
for $n\equiv 2\mod 3$
$$
\f_{\r_2}(m,k,n)=\Tr\left[
\Matt 0{e^{\pi i m}}0{e^{\pi i k}}0000{e^{\pi i(k+m)}}
\Matt 010001100
\right]=
$$
$$
=
\Tr\Matt 00{e^{\pi i m}}0{e^{\pi i k}}0{e^{\pi i(k+m)}}00
=e^{\pi i k}=\left\{
\begin{array}{rl}
1,& \mbox{ if } k \mbox{ is even }\\
-1,& \mbox{ if } k \mbox{ is odd }
\end{array}
\right. .
$$
$\f_{\r_2}$ is $3$-periodical in $n$, while the characteristic functions
of $B_i$ are $6$-periodical. For $j=0,\dots,5$ one has
$$
\begin{array}{llll}
\f_{\r_2}|_{B_1\cap L_0}\equiv 1, &
\f_{\r_2}|_{B_2\cap L_0}\equiv -1, &
\f_{\r_2}|_{B_3\cap L_1}\equiv 1, &
\f_{\r_2}|_{B_4\cap L_1}\equiv -1, \\
\f_{\r_2}|_{B_1\cap L_2}\equiv 1, &
\f_{\r_2}|_{B_2\cap L_2}\equiv -1, &
\f_{\r_2}|_{B_3\cap L_3}\equiv 1, &
\f_{\r_2}|_{B_4\cap L_3}\equiv -1, \\
\f_{\r_2}|_{B_1\cap L_4}\equiv 1, &
\f_{\r_2}|_{B_2\cap L_4}\equiv -1, &
\f_{\r_2}|_{B_3\cap L_5}\equiv 1, &
\f_{\r_2}|_{B_4\cap L_5}\equiv -1,
\end{array}
$$
so $\f_{\r_2}|_{B_1\cup B_3}\equiv 1$, $\f_{\r_2}|_{B_2\cup B_4}\equiv -1$.
the determinant of the values of the functions $\f_{\r_1}$, $\f_\pi$,
$\f_{\r_2}$, $\f_{\r_2\otimes\pi}$ on $B_1$, $B_2$, $B_3$, $B_4$ is
$$
\det \left(
\begin{array}{rrrr}
1 & 1 & 1 & 1 \\
1 & 1 & -1 & -1 \\
1 & -1 & 1 & -1 \\
1 & -1 & -1 & 1
\end{array}
\right)
=-8\ne 0.
$$
Hence, they are linearly independent.

Nevertheless, there are infinite-dimensional irreducible
$\wh\phi$-invariant representations. E.g. we have a representation
$\r$ of $G$ on $L^2(\T^2)$ with the respect to the Lebesgue measure,
with $\r(m,k,0)$ be the multiplier by characters in the appropriate
points and $\r(0,0,1)$ is $\a$, in the same manner as for $\r_2$.
This representation can
be also obtained by inducing from the trivial representation of the subgroup
$\{\alpha^n\}\subset G$ and
Fourier transform. The irreducibility of this induced representation can also
be proved using \cite{Saito}.

This disproves the conjecture of Fel'shtyn and Hill \cite{FelHill},
who supposed
that the Reidemeister number equals to the number of fixed points of
$\wh \phi$ on $\wh G$

This representation is not traceable, but
one can nevertheless try to calculate (\ref{eq:trsg}).
We will do this in some fixed base, because in general operators
are not of trace class.
Let us choose
an orthonormal base of $L^2(\T^2)$ formed by
$\e_{st}(x,y)=e^{2\pi i(sx+ty)}$, $x,y\in [0,1]$. The intertwining
operator $S$ is generated by $\mu$. Then
$$
\Tr(S\r(m,k,n))=\sum_{s,t} \int_0^1\int_0^1 (\r(m,k,0)
\r(0,0,n) \e_{s,t})(\mu(x,y))
\ov{\e_{s,t}}(x,y) \,dx\,dy=
$$
$$
=\sum_{s,t} \int_0^1\int_0^1 e^{2\pi i \<(m,k),\mu(x,y)\>}
(\r(0,0,n) \e_{s,t})(\mu(x,y))
\e_{-s,-t}(x,y) \,dx\,dy=
$$
$$
=\sum_{s,t} \int_0^1\int_0^1 e^{2\pi i \<(m,k),\mu(x,y)\>}
\e_{s,t}(\a^{n}\mu(x,y))
\e_{-s,-t}(x,y) \,dx\,dy=
$$
$$
=\sum_{s,t} \int_0^1\int_0^1
e^{2\pi i \<(m,k)+(\a^{n}-\mu)(s,t),\mu(x,y)\>} \,dx\,dy=
\sum_{(m,k)=(\mu-\a^n)(s,t)} 1.
$$
For $n=0$ this equals $1$ for $m+k$ even and $0$ for $m+k$ odd.
Hence, $\f_\r|_{B_1\cap L_0}=1$, $\f_\r|_{B_2\cap L_0}=0$. For $n=1$
the equality takes the form $(m,k)=(-2s-2t,-t)$. Hence,
$\f_\r|_{B_3\cap L_1}=1$, $\f_\r|_{B_4\cap L_1}=0$.
If $n=2r$, the equality takes form
$$
(m,k)=(\mu-\a^{2r})(s,t),\quad
(m,k)=(\a^r \mu \a^r -\a^{2r})(s,t),
$$
$$
\a^{-r} (m,k)=(\mu-\a^0) \a^{r}(s,t).
$$
Since $\a$ is an automorphism of $\Z\oplus\Z$ and by the description
of $B_i$ via the action of $\a$, we obtain that
$\f_\r|_{B_1}=1$, $\f_\r|_{B_2}=0$. Similarly, for $n$ odd.
So, $\f_\r$ is well defined and
$$
\f_\r=\chi_{B_1}+\chi_{B_3}=\frac 12 (\f_{\r_1}+\f_{\r_2}).
$$

Of course, there are also $\wh\phi$-invariant traceable factor
representations of this group $G$, e.g. the regular representation.
Since its kernel is trivial, evidently all twisted-invariant
functionals can be pulled back from it.

\section{Discussion, prospectives, conjectures}

This example shows that the conjecture of \cite{FelHill} for general
groups can survive only after eliminating badly separated points in
$\wh G$. A partial solution of the problem along this direction is
obtained in Theorem \ref{teo:weakburn}.
As it is known, the unitary
dual of a group of polynomial growth is well separated, this gives hope to
get an analytical proof of a version of twisted Burnside theorem for these groups.
Also, while the present paper was under refereeing A.F. and E.T. have obtained
 the following general theorem:
for any automorphism of any almost polycyclic group $R(\phi)=S_f(\phi)$ if one of
them is finite, where $S_f(\phi)$ is the number of fixed points of $\wh\phi$ on the
finite-dimensional part of $\wh G$\
(\cite{polyc}, see also the survey \cite{FelTroObzo} for a number of related
statements, counterexamples, etc.).


\begin{thebibliography}{10}

\bibitem{Arthur}
J.~Arthur and L.~Clozel, \emph{Simple algebras, base change, and the advanced
  theory of the trace formula}, Princeton University Press, Princeton, NJ,
  1989. \MR{90m:22041}

\bibitem{baron}
A.~O. Barut and R.~R{\polhk{a}}czka, \emph{Theory of group representations and
  applications}, second ed., World Scientific Publishing Co., Singapore, 1986.
  \MR{88c:22013}

\bibitem{DaunsHofmann}
J.~Dauns and K.~H. Hofmann, \emph{Representations of rings by continuous
  sections}, Mem. Amer. Math. Soc., vol.~83, Amer. Math. Soc., Providence, RI,
  1968.

\bibitem{FelshB}
A.~Fel'shtyn, \emph{Dynamical zeta functions, {N}ielsen theory and
  {R}eidemeister torsion}, Mem. Amer. Math. Soc. \textbf{147} (2000), no.~699,
  xii+146. \MR{2001a:37031}

\bibitem{FelGon}
A.~Fel'shtyn and D.~Gon{\c{c}}alves, \emph{Reidemeister numbers of
  {B}aumslag-{S}olitar groups}, E-print arXiv:math.GR/0405590, 2004, (to appear
  in \emph{Algebra and Discrete Mathematics}).

\bibitem{FelHill}
A.~Fel'shtyn and R.~Hill, \emph{The {R}eidemeister zeta function with
  applications to {N}ielsen theory and a connection with {R}eidemeister
  torsion}, $K$-Theory \textbf{8} (1994), no.~4, 367--393. \MR{95h:57025}

\bibitem{FelBanach}
\bysame, \emph{Dynamical zeta functions, congruences in {N}ielsen theory and
  {R}eidemeister torsion}, Nielsen theory and Reidemeister torsion (Warsaw,
  1996), Polish Acad. Sci., Warsaw, 1999, pp.~77--116. \MR{2001h:37047}

\bibitem{FelTroObzo}
A.~Fel'shtyn and E.~Troitsky, \emph{Geometry of {R}eidemeister classes and
  twisted {B}urnside theorem}, (to appear in \emph{K-Theory}).

\bibitem{polyc}
\bysame, \emph{Twisted {B}urnside theorem}, Preprint 46, Max-Planck-Institut
  f{\"u}r Mathematik, 2005.

\bibitem{FelTro}
\bysame, \emph{A twisted {B}urnside theorem for countable groups and
  {R}eidemeister numbers}, Noncommutative Geometry and Number
  Theory (C.~Consani, M.~Marcolli, eds.), Vieweg,
  Braunschweig, 2006, pp.~141--154, (Preprint MPIM2004-65, math.RT/0606155).

\bibitem{FelPOMI}
A.~L. Fel'shtyn, \emph{The {R}eidemeister number of any automorphism of a
  {G}romov hyperbolic group is infinite}, Zap. Nauchn. Sem. S.-Peterburg.
  Otdel. Mat. Inst. Steklov. (POMI) \textbf{279} (2001), no.~6 (Geom. i
  Topol.), 229--240, 250. \MR{2002e:20081}

\bibitem{gowon}
D.~Gon{\c{c}}alves and P.~Wong, \emph{Twisted conjugacy classes in exponential
  growth groups}, Bull. London Math. Soc. \textbf{35} (2003), no.~2, 261--268.
  \MR{2003j:20054}

\bibitem{Groth}
A.~Grothendieck, \emph{Formules de {N}ielsen-{W}ecken et de {L}efschetz en
  g\'eom\'etrie alg\'ebrique}, S\'eminaire de G\'eom\'etrie Alg\'ebrique du
  {B}ois-{M}arie 1965-66. {SGA} 5, Lecture Notes in Math., vol. 569,
  Springer-Verlag, Berlin, 1977, pp.~407--441.

\bibitem{Jiang}
B.~Jiang, \emph{Lectures on {N}ielsen fixed point theory}, Contemp. Math.,
  vol.~14, Amer. Math. Soc., Providence, RI, 1983.

\bibitem{Kirillov}
A.~A. Kirillov, \emph{Elements of the theory of representations},
  Springer-Verlag, Berlin Heidelberg New York, 1976.

\bibitem{ll}
G.~Levitt and M.~Lustig, \emph{{Most automorphisms of a hyperbolic group have
  very simple dynamics.}}, Ann. Scient. \'Ec. Norm. Sup. \textbf{33} (2000),
  507--517.

\bibitem{Saito}
M.~Saito, \emph{Repr\'esentations unitaires monomiales d'un groupe discret, en
  particulier du groupe modulaire}, J. Math. Soc. Japan \textbf{26} (1974),
  no.~3, 464--482.

\bibitem{Shokra}
S.~Shokranian, \emph{The {S}elberg-{A}rthur trace formula}, LNM No. 1503,
Springer-Verlag,
  Berlin, 1992, Based on lectures by James Arthur. \MR{93j:11029}

\bibitem{ncrmkwb}
E.~Troitsky, \emph{Noncommutative {R}iesz theorem and weak
  {B}urnside type theorem on twisted conjugacy},
Funct. Anal. Appl. \textbf{40} (2006), No. 2, 44--54. In Russian,
English translation to appear in \emph{Funct. Anal. Appl.}
  (Preprint 86, Max-Planck-Institut f{\"u}r
  Mathematik, 2004).

\end{thebibliography}

\def\cprime{$'$} \def\dbar{\leavevmode\hbox to 0pt{\hskip.2ex \accent"16\hss}d}
  \def\polhk#1{\setbox0=\hbox{#1}{\ooalign{\hidewidth
  \lower1.5ex\hbox{`}\hidewidth\crcr\unhbox0}}}
\providecommand{\bysame}{\leavevmode\hbox to3em{\hrulefill}\thinspace}
\providecommand{\MR}{\relax\ifhmode\unskip\space\fi MR }
\providecommand{\MRhref}[2]{%
  \href{http://www.ams.org/mathscinet-getitem?mr=#1}{#2}
}
\providecommand{\href}[2]{#2}

\end{document}